\documentclass[%
 reprint,
 amsmath,amssymb,
 aps,
]{revtex4-2}

\usepackage{graphicx}
\usepackage{dcolumn}
\usepackage{bm}
\usepackage{color}
\usepackage{esint}

\begin{document}

\preprint{APS/123-QED}

\title{Sea urchin sperm exploit extremum seeking control to find the egg}
\author{Mahmoud Abdelgalil}
\email{maabdelg@uci.edu}
\affiliation{Mechanical and Aerospace Engineering Department,  University of California at Irvine.}
\author{Yasser Aboelkassem}
\email{yassera@umich.edu}
\affiliation{College of Innovation and Technology, University of Michigan at Flint, Flint MI 48502 \\ Michigan Institute for Data Science, University of Michigan, Ann Arbor MI 48109}
\author{Haithem Taha}
\email{hetaha@uci.edu}
\affiliation{Mechanical and Aerospace Engineering Department, University of California at Irvine.
}
\date{\today}

\begin{abstract}
Sperm cells perform extremely demanding tasks with minimal capabilities. The cells must quickly navigate in a noisy environment to find an egg within a short time window for successful fertilization without any global positioning information. Many research efforts have been dedicated to derive mathematical principles that explain their superb navigation strategy. Here we show that the navigation strategy of sea urchin sperm, also known as helical klinotaxis, is a natural implementation of a well-established adaptive control paradigm known as extremum seeking. This bridge between control theory and the biology of taxis in microorganisms is expected to deepen our understanding of the process. For example, the formulation leads to a coarse-grained model of the signaling pathway that offers new insights on the peculiar switching-like behavior between high and low gain steering modes observed in sea urchin sperm. Moreover, it may guide engineers in developing bio-inspired miniaturized robots with minimal sensors.
\end{abstract}

\maketitle
\section{Introduction}
Source seeking, a well-studied topic in the control community \cite{krstic2008extremum}, is the problem of locating an object that emits a scalar measurable signal (e.g. chemical concentration, sound, heat, etc.), typically without global positioning information. Many organisms are routinely faced with the source seeking problem. A well studied example is that of sperm chemotaxis \cite{friedrich2007chemotaxis,alvarez2014computational}. To locate an egg in open water, sea urchin sperm evolved to swim up the gradient of the concentration field established by the diffusion of a species-specific chemoattractant, a sperm-activating peptide (SAP), secreted by the eggs \cite{alvarez2014computational}. Unlike the inherently stochastic bacterial chemotaxis, the navigation strategy of sea urchin sperm can be reasonably described in a deterministic fashion; the cells employ the mean curvature of the flagellum, regulated by intracellular calcium, as a steering mechanism to swim in circular paths that drift in the direction of the gradient in 2D, and in helical paths that align with the gradient in 3D \cite{crenshaw1993orientation,friedrich2007chemotaxis,friedrich2009steering,jikeli2015sperm}. This feedback mechanism is mediated by a complex signaling pathway that regulates the influx and efflux of calcium in the cell \cite{kaupp2003signal,priego2020modular}. 

In this letter, we revisit sperm chemotaxis from the perspective of control theory. We frame the search for the egg as a source seeking problem, then we show that the 3D navigation strategy of sea urchin sperm, also known as helical klinotaxis, is in fact a natural implementation of a well established adaptive control paradigm known as extremum seeking \cite{krstic2008extremum,scheinker2013model,scheinker2014extremum}. We illustrate this novel connection by establishing a one-to-one correspondence between the key components of the navigation strategy of sea urchin sperm cells and the hallmark features of an extremum seeking solution to the source seeking problem.
Based on this formulation, we propose a coarse-grained minimal dynamical description that captures the crucial features of the chemotactic signaling pathway, including the peculiar behavior of sea urchin sperm cells where they seem to switch between two distinct navigation modes: i) the `on-response' which is a low-gain steering mode when the average velocity vector of the cell is mostly aligned with the gradient, and ii) the `off-response' which is a high-gain steering mode otherwise \cite{jikeli2015sperm}. The proposed description improves upon previous models \cite{friedrich2007chemotaxis,jikeli2015sperm,kromer2018decision} as it does not employ an explicit discontinuous switching logic to explain the switching-like phenomenon. Instead, the behavior naturally arises as a consequence of the motion pattern and a time-scale separation between the proposed dynamics of the signaling pathway and the average motion. In particular, the proposed model does not exploit any information other than the perceived instantaneous local concentration which is readily available to the cells through the SAP receptors located on the flagellum.
\section{A Primer on Extremum Seeking}
\begin{figure}[t]
    \centering
        \includegraphics[width=0.9\columnwidth]{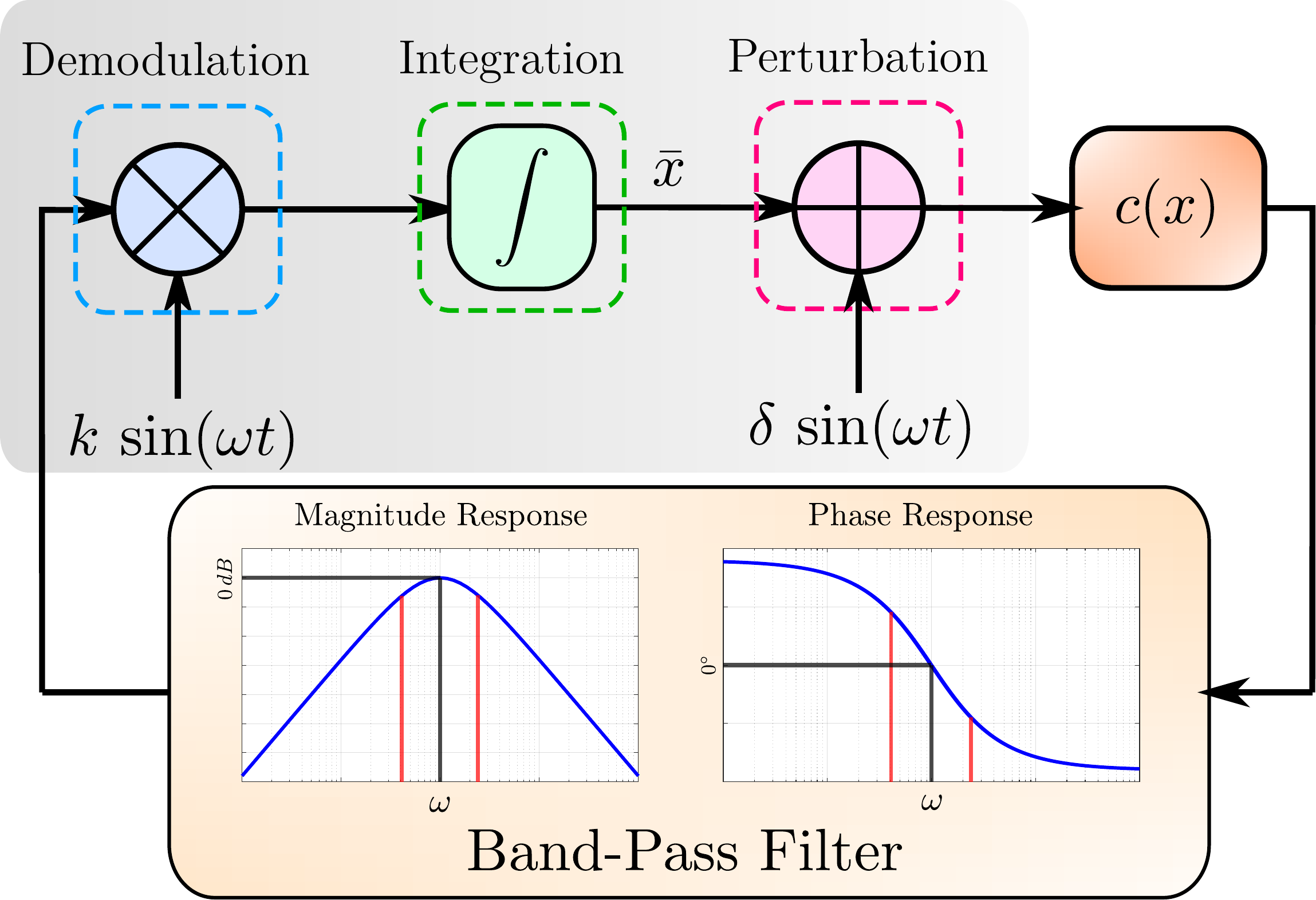}
        \caption{A block diagram description of the simplest extremum seeking control scheme as represented by equations (\ref{eq:simple_es}).} \label{fig:blk_diag_es}
\end{figure}

{\color{black} We begin with a brief exposition of Extremum Seeking (ES) control. ES is an adaptive control technique designed to steer a dynamical system towards the extremum of an objective function that depends on the state of the system, without access to information about the gradient of the function (only the value of the objective function is available for measurement at each instant in time). The first ES control law can be traced back to the century old paper due to Leblanc \cite{leblanc1922electrification}, but the recent interest in ES control was sparked by Krsti\'c's seminal paper \cite{krstic2000stability}. In the simplest setting, an ES controller is designed to find the optimal value of a single-variable static objective function by dynamically estimating the gradient. Let $c(x)$ be the objective function, and consider the following dynamical system \cite{ariyur2003real}:
\begin{subequations}\label{eq:simple_es}
    \begin{gather}\label{eq:simple_es_1}
        x = \bar{x} + \delta\,\sin(\omega t), \qquad
        \dot{\bar{x}}= 2\,\zeta\,k\,\sin(\omega t), \\
        \label{eq:simple_es_2}
        \dot{\zeta}_1= \omega\,(\zeta_2 - \zeta_1), \; \dot{\zeta}_2= \omega\,(c(x) - \zeta_2), \; \zeta = \zeta_2-\zeta_1,
    \end{gather}
\end{subequations}
which is depicted in the block diagram \cite{abramovici2000feedback} presented in Fig.\ref{fig:blk_diag_es}, where $\bar{x}$ is the estimate of the optimal value of the independent variable $x$, $\zeta_1$ and $\zeta_2$ are the states of a band-pass filter centered around $\omega$, and $k$, $\omega$, $\delta$ are constants. The flow of the block diagram in Fig.\ref{fig:blk_diag_es} can be traced as follows. First, a sinusoidal perturbation is injected to sample the objective function near the estimate $\bar{x}$:
\begin{align}\label{eq:cost_approx}
    c(x)&= c(\bar{x}) + \frac{dc(\bar{x})}{d\bar{x}}\,\delta \,\sin(\omega t) + O(\delta^2).
\end{align}
We observe how the gradient appears as the amplitude of the sinusoidal perturbation. In engineering terms, injecting the perturbation around the current estimate $\bar{x}$ `modulates' the local gradient information on the amplitude of the sinusoidal `carrier' signal $\sin(\omega t)$. Therefore, to extract the sinusoidal signal that carries the gradient information, the measured objective function $c(x)$ goes through a band-pass filter centered around the frequency $\omega$ as defined by equation (\ref{eq:simple_es_2}). The output of the filter $\zeta$ can be approximated in the quasi-steady sense by:
\begin{align}
    2\zeta \approx 2\zeta_{\text{QS}} = \frac{dc(\bar{x})}{d\bar{x}}\,\delta \,\sin(\omega t),
\end{align}
Next, the gradient information is `demodulated' (i.e., extracted from the carrier signal) through multiplication with a sinusoidal signal having the same frequency and phase as the carrier signal:
\begin{align}\label{eq:demodulated}
    2\,\zeta\,k\sin(\omega t)&= \frac{dc(\bar{x})}{d\bar{x}}\,\delta\,k - \frac{dc(\bar{x})}{d\bar{x}}\,\delta\,k\,\cos(2\omega t),
\end{align}
where the time-average of the right hand side of equation (\ref{eq:demodulated}) is non-zero. Finally, the demodulated gradient information is used in adjusting the current estimate $\bar{x}$. Through a simple averaging argument, we obtain that the estimate $\bar{x}$ on average evolves according to:
\begin{align}
    \dot{\bar{x}} &\approx \overline{2\,\zeta_{\text{QS}}\,k\,\sin(\omega\,t)} = k\,\delta\,\frac{d c(\bar{x})}{d\bar{x}},
\end{align}
where the overline $\overline{\hphantom{\,\,\,}\vphantom{h}}$ indicates the time average of the overlined quantity. That is, the estimate $\bar{x}$ evolves, in a quasi-steady average sense, along the gradient of the objective function under the extremum seeking control law (\ref{eq:simple_es}). The interested reader is referred to \cite{krstic2000stability,ariyur2003real,tan2010extremum} for more details. }
\section{Modeling the sperm motion} 
We now turn our attention to the motion of the sperm cell. Swimming in a low Reynolds number is dominated by viscous forces, which enables the use of kinematic models as a good approximation to the motion of micro-swimmers, including sperm cells \cite{friedrich2010high}. The kinematics of a rigid body are given by:
\begin{align}
    \label{eq:orig_kin} \dot{\textbf{p}} &= \textbf{R}\textbf{v}, & \dot{\textbf{R}} &= \textbf{R}\bm{\widehat{\Omega}},
\end{align}
where the vectors $\textbf{v}$ and $\bm{\Omega}$ are the linear and angular velocity vectors in the body frame, $\bm{\widehat{\Omega}}$ denotes the skew-symmetric matrix corresponding to the angular velocity vector $\bm{\Omega}$, $\textbf{p}$ is the instantaneous position of the body with respect to the origin of a fixed frame of reference, and $\textbf{R}$ is the instantaneous rotation matrix that relates the body frame to the fixed frame. In sea urchin sperm, the mean curvature and torsion of the flagellar beating pattern, which are regulated by the chemotactic signaling pathway, control the angular velocities in the body frame \cite{crenshaw1993orientation,alvarez2014computational}. A common model of the effect of the chemotactic signaling pathway on the swimming kinematics of sea urchin sperm is given by the relations:
\begin{align}\label{eq:sperm_kin_model}
    \textbf{v} &= \left[\begin{array}{ccc} v & 0 & 0 \end{array}\right]^\intercal, &\bm{\Omega} &= \left[\begin{array}{ccc} \omega_\parallel & 0 & \omega_\perp \end{array}\right]^\intercal,
\end{align}
where $v>0$ is constant, and the angular velocity components $\omega_\parallel$ and $\omega_\perp$ are given by:
\begin{align}\label{eq:omega_def}
    \omega_\parallel &= \omega_{\parallel 0} + \omega_{\parallel 1}  \eta , & \omega_\perp &= \omega_{\perp 0} + \omega_{\perp 1}  \eta ,
\end{align}
with $\omega_{\perp 0},\omega_{\perp 1},\omega_{\parallel 0},\omega_{\parallel 1}$ as constant coefficients, and $\eta$ is a dynamic feedback term regulated by the signaling pathway \cite{crenshaw1996new,friedrich2007chemotaxis,jikeli2015sperm}.

\section{An extremum seeking loop} 
The constant forward velocity $v>0$, along with the constant angular velocity components $\omega_{\parallel 0}$ and $\omega_{\perp 0}$, lead to a periodic swimming pattern, a helical trajectory, which injects periodic perturbations into the instantaneous position and orientation of the cell. The sign of $\omega_{\parallel 0}$ and $\omega_{\perp 0}$ determine the handedness of the helical trajectory. For simplicity, we consider the case in which both $\omega_{\parallel 0}$ and $\omega_{\perp 0}$ are positive. {\color{black} We define the average instantaneous position $\bar{\textbf{p}}$ and orientation $\bar{\textbf{R}}$ of the cell as:
\begin{subequations}\label{eq:avg_mot}
    \begin{gather}
        \begin{aligned}
            \textbf{R}_0(t) &= \text{exp}\big(\bm{\widehat{\Omega}}_0 t\big), & \bar{\textbf{R}}&= \textbf{R} \textbf{R}_0(t)^\intercal,
        \end{aligned}\\
        \bar{\textbf{p}}= \textbf{p}- \bar{\textbf{R}}\bm{\delta}(t),
    \end{gather}
\end{subequations}
where the vector $\bm{\Omega}_0$ is given by:
\begin{align}
    \bm{\Omega}_0 &= \left[\begin{array}{ccc} \omega_{\parallel 0} & 0 & \omega_{\perp 0} \end{array}\right]^\intercal
\end{align}
and the time-periodic vector $\bar{\textbf{R}} \bm{\delta}(t)$ is the perturbation in the position due to the helical swimming pattern, and is defined by:
\begin{align}
    \textbf{v}_m&= \overline{\textbf{R}_0(t)}\textbf{v}, & \bm{\delta}(t) &= \int \left(\textbf{R}_0(t)\textbf{v} - \textbf{v}_m\right)dt.
\end{align}
Direct computations show that:
\begin{align}
    \lvert \textbf{v}_m\lvert= v\omega_{\parallel 0}/\omega,\;  \lvert\bm{\delta}(t)\lvert = v\omega_{\perp 0}/\omega^2,\; \textbf{v}_m^\intercal \bm{\delta}(t) = 0,
\end{align}
where $\omega = \lvert \bm{\Omega}_0\lvert$ is the frequency of the periodic motion. In particular, the periodic perturbation $\bar{\textbf{R}} \bm{\delta}(t)$ and the direction of the average motion $\bar{\textbf{R}} \textbf{v}_m$ are orthogonal. The evolution of the average motion variables $\bar{\textbf{p}}$ and $\bar{\textbf{R}}$ is governed by the following system of differential equations with periodic coefficients:
\begin{subequations}\label{eq:kin}
\begin{align}
\label{eq:rot_kin}
    \dot{\bar{\textbf{R}}} &=\bar{\textbf{R}}\widehat{\bm{\Omega}}_\eta(t)\eta, & \bm{\Omega}_\eta(t)&= \textbf{R}_0(t){\bm{\Omega}}_1, \\
\label{eq:pos_kin}
    \dot{\bar{\textbf{p}}}&= \bar{\textbf{R}}\textbf{v}_\eta(t)\eta +  \bar{\textbf{R}}\textbf{v}_m, & \textbf{v}_\eta(t)&=\bm{\delta}(t)\times\bm{\Omega}_\eta(t),
\end{align} 
\end{subequations}
where the vector $\bm{\Omega}_1$ is given by:
\begin{align}
    \bm{\Omega}_1 &= \left[\begin{array}{ccc} \omega_{\parallel 1} & 0 & \omega_{\perp 1} \end{array}\right]^\intercal.
\end{align}
In the absence of feedback  (i.e. when $\eta = 0$), the vector $\bar{\textbf{R}}\textbf{v}_m$ is the average velocity vector of the cell, and its direction is along the axis of the helical trajectory. As such, the instantaneous local SAP concentration $c(\textbf{p})$, which stimulates the receptors along the flagellum \cite{alvarez2014computational} and is assumed here to be a smooth function of position, can be approximated by its first order Taylor series, over a short time duration, in terms of the average motion variables $\bar{\textbf{p}}$ and $\bar{\textbf{R}}$ using the relations (\ref{eq:avg_mot}) and (\ref{eq:pos_kin}):
\begin{equation}\label{eq:stim}
    \begin{aligned}
        c(\textbf{p}) \approx c(\bar{\textbf{p}}_0) + \nabla c(\bar{\textbf{p}})^\intercal\bar{\textbf{R}}\textbf{v}_m \Delta t+ \nabla c(\bar{\textbf{p}})^\intercal\bar{\textbf{R}}\bm{\delta}(t),
    \end{aligned}
\end{equation}
where $\bar{\textbf{p}}_0$ is the average position of the cell at the initial time, and we assume that $\lvert \bm{\delta}(t)\lvert \ll 1$. Clearly, the injection of periodic perturbations due to the helical swimming pattern modulates the local gradient information on the amplitude of the periodic perturbations $\bar{\textbf{R}}\bm{\delta}(t)$. In other words, the periodic motion pattern acts as a carrier signal upon which the gradient information is modulated, in an identical manner to the perturbation stage of the ES control loop.\\

Since the perturbation vector $\bar{\textbf{R}}\bm{\delta}(t)$ and the average velocity vector $\bar{\textbf{R}}\textbf{v}_m$ are orthogonal, the amplitude of the periodic signal $\nabla c(\bar{\textbf{p}})^\intercal\bar{\textbf{R}}\bm{\delta}(t)$ is proportional to the orthogonal misalignment between the average direction of motion of the cell as defined by $\bar{\textbf{R}}\textbf{v}_m$ and the local gradient. Hence, the signaling pathway is ought to produce a feedback signal $\eta$ that eliminates this misalignment. It is well known that microorganisms that swim in helical trajectories, including sea urchin sperm, can align the axis of their helical trajectory with the gradient by periodically varying the angular velocities of the cell with the same frequency of the helical trajectory \cite{crenshaw1993orientation,crenshaw1996new}. That is, a sperm cell can align its average direction of motion with the gradient provided that the signaling pathway is able to extract the periodic component of the instantaneous local concentration. For successful chemotaxis, the signaling pathway is ought to play a similar role of the filter in the ES control loop, by extracting the periodic component in the instantaneous concentration with the same frequency of the periodic perturbation $\bar{\textbf{R}}\bm{\delta}(t)$ that carries the gradient information. This implication about the nature of the dynamics of the signaling pathway is one of the main outcomes of the connection between chemotaxis and ES, as proposed in this paper. \\

Going back to the governing equations of the average kinematics (\ref{eq:kin}), we see that the feedback signal $\eta$ multiplies the periodic feedback coefficients $\textbf{v}_\eta(t)$ and $\bm{\Omega}_\eta(t)$. Consequently, the local gradient information carried on the periodic component in the signal $\eta$ is `demodulated' into the non-zero average component of the product signals $\bm{\Omega}_\eta(t)\,\eta$ and $\textbf{v}_\eta(t)\,\eta$, similar to the demodulation stage of the ES control loop.

Finally, the demodulated local gradient information passes through the kinematics of the motion represented by equations (\ref{eq:kin}), which is responsible for biasing the motion in the direction of the gradient. The closed-loop behavior of the nonholonomic integrator defined by the kinematics is investigated in the next section. A block diagram description of the dynamical equations (\ref{eq:kin}) representing the navigation strategy of sea urchin sperm is shown in Fig.\ref{fig:blk_diag_ct}, where the special integration symbol $\fint$ denotes the nonholonomic kinematic integrator corresponding to the equations (\ref{eq:kin}). The isomorphism between the block diagrams in Fig.\ref{fig:blk_diag_es} and Fig.\ref{fig:blk_diag_ct} clearly reveals the connection between sperm chemotaxis and extremum seeking. }\\

It is worth mentioning that the 2D version of the model (\ref{eq:orig_kin})-(\ref{eq:omega_def}) (i.e. when $\omega_\parallel = 0$ and the motion is restricted to a plane) is a well-studied kinematic model in the control community known as the unicycle model. Remarkably, the trajectories generated by an ES-based algorithm for the unicycle model, which was recently proposed in \cite{scheinker2013model,scheinker2014extremum} independently from the literature on sperm chemotaxis, are astonishingly similar to the actual trajectories of sea urchin sperm in shallow observation chambers \cite{bohmer2005ca2+}.
\begin{figure}[t]
    \centering
        \includegraphics[width=\columnwidth]{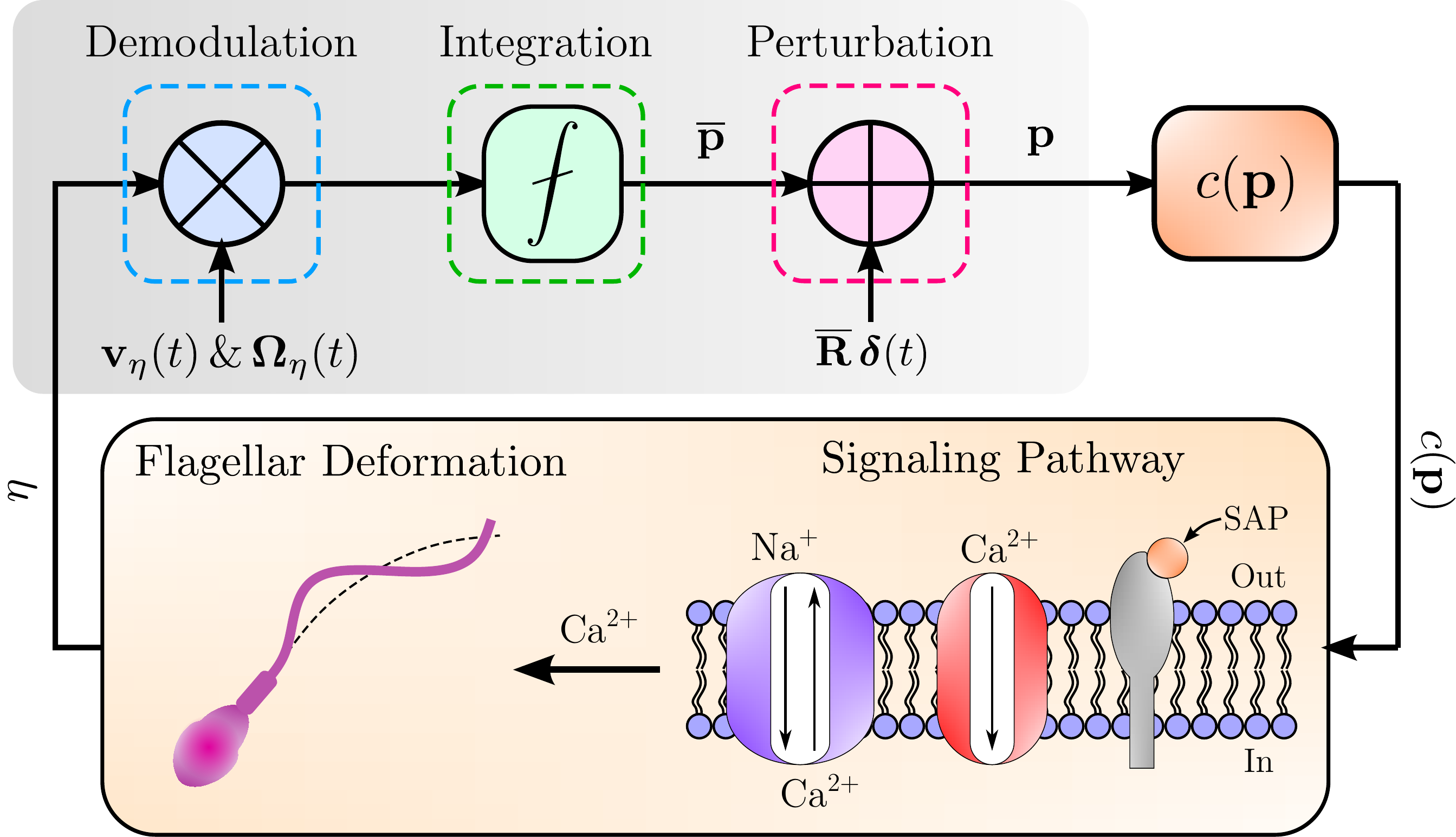}
        \caption{A block diagram description of equations (\ref{eq:avg_mot}) and (\ref{eq:kin}). The swimming pattern injects periodic perturbations into the instantaneous position of the cell which leads to oscillations in the instantaneous local SAP concentration. The signaling pathway relays these periodic perturbations to the angular velocities through flagellar deformation. Then, the periodic feedback coefficients (i.e. $\textbf{v}_\eta$ and $\bm{\Omega}_\eta$) of the dynamics of average motion demodulate the gradient information carried by the feedback signal $\eta$ through signal multiplication. Finally, the kinematic integrator $\fint$ biases the motion in the direction of the gradient.} \label{fig:blk_diag_ct}
\end{figure}

\begin{figure}[t]
    \centering
    \includegraphics[width=\columnwidth]{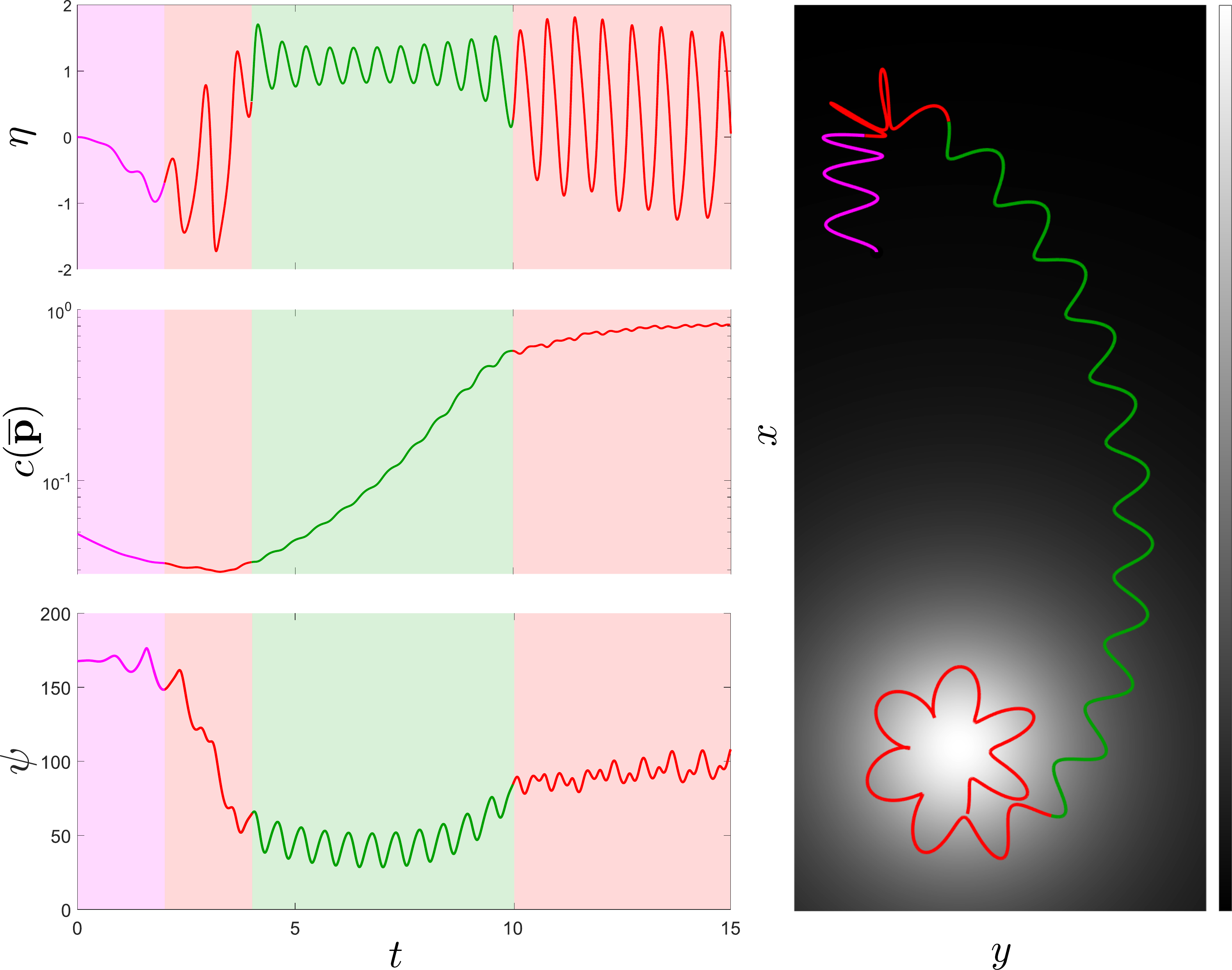}\vspace*{0.1in}
    \begin{tabular}{|c|c |c |c |c|c|c|c|c|c|} 
     \hline
     $v$ & $\omega_{\parallel 0}$ & $\omega_{\parallel 1}$ & $\omega_{\perp 0}$ & $\omega_{\perp 1}$ & $\mu$ & $\sigma$ & $\eta(t_0)$ & $\zeta(t_0)$ & $\rho(t_0)$ \\
     \hline
    3.07 & 3.07 & 2.30 & 8.91 & 1.00 & 3.07 & 9.42 & 0.00 & $c(\textbf{p}_0)$ & 5.00 \\
     \hline
    \end{tabular}
    \caption{The three cases of the behavior of the signaling pathway illustrated on a sample trajectory projected on the $xy$-plane, the response $\eta$, the average instantaneous concentration $c(\bar{\textbf{p}})$, and the angle $\psi = \cos^{-1}(\bar{\textbf{h}}^\intercal\check{\nabla}c)$ (in degrees) between the gradient and the average direction of motion $\bar{\textbf{h}}$, in a radial concentration field $c(\textbf{p})=1/(1+0.5\lvert \textbf{p}\lvert^2)$. The initial position is taken as $\textbf{p}_0=(6,1,0)$, and the initial orientation is $\textbf{R}(t_0)=\text{exp}(2\pi\hat{\textbf{e}}_2/5)$, where $\hat{\textbf{e}}_2 = (0,1,0)$. The rest of the initial conditions and parameter values are in the table. 
    }

    \label{fig:sample_traj}
\end{figure}
\section{Chemotactic Response and Closed Loop Behavior}
The evident one-to-one correspondence between the key components of the navigation strategy of sea urchin sperm and ES control immediately clarifies the role of the signaling pathway: it must act as an adaptive band-pass filter attuned to the frequency of the swimming pattern of the cell. Motivated by this observation, and building upon previous phenomenological models \cite{friedrich2007chemotaxis,jikeli2015sperm}, we propose the following coarse-grained dynamical description of the signaling pathway: 
\begin{subequations}
    \label{eq:pthway}
    \begin{align}
        \label{eq:pthway_1}
        \sigma\,\dot{\xi}&= s(t) - \xi,\\
        \label{eq:pthway_2}
        \mu\,\dot{\eta}&= \rho\,\dot{\xi}-\eta^3,\\
        \label{eq:pthway_3}
        \mu\,\dot{\rho} &= \rho-\rho \,\eta^2,
    \end{align} 
\end{subequations}
where $\mu$ and $\sigma$ are positive constants such that $\sigma < \mu$, and $s(t)$ is the input to the model, which represents the time-varying external stimulus to which the pathway is exposed due to the binding of SAP molecules with the receptors. Without accounting for noise, the stimulus $s(t)$ is customarily approximated by:
\begin{align}
    s(t) \approx \lambda \, c(\textbf{p})
\end{align}
for some positive proportionality constant $\lambda$ \cite{friedrich2007chemotaxis,jikeli2015sperm}. The proposed model possesses three essential dynamical features: excitation, relaxation, and adaptation. The excitation is modelled by equation (\ref{eq:pthway_1}), which acts as a differentiator that detects changes in the local concentration. The relaxation is modelled by equation (\ref{eq:pthway_2}), which brings the response $\eta$ back to resting levels when there is no change in the stimulus. Finally, the adaptation is modelled by equation (\ref{eq:pthway_3}), which adjusts the sensitivity of the pathway to the stimulus. A sample trajectory of the equations (\ref{eq:orig_kin})-(\ref{eq:omega_def}) and (\ref{eq:pthway}) is shown in Fig.\ref{fig:sample_traj} along with time-history of the average local concentration $c(\bar{\textbf{p}})$, the steering response $\eta$ and the angle $\psi$ between the gradient and the direction of motion $\bar{\textbf{h}}$. \\

{\color{black}We now analyze the closed loop behavior when the dynamics of the pathway is given by the proposed dynamical system (\ref{eq:pthway}). The details of calculations in this section can be found in the SI appendix.  In the parametric regime where $\sigma\lvert \textbf{v}_m\lvert\ll \mu\lvert \textbf{v}_m\lvert\ll\sigma\lvert \bm{\Omega}_0\lvert \approx O(1)$, there is a large time-scale separation between the dynamics of average motion (\ref{eq:kin}) in the absence of feedback and the dynamics of the pathway (\ref{eq:pthway}). Consequently, we may approximate the response $\eta$ due to the time-varying local concentration (\ref{eq:stim}) by the quasi-steady response:
\begin{align}\label{eq:eta_qs}
   \eta_{\text{QS}} = \frac{  \check{\nabla}_{\hspace*{-0.02in}\parallel}c + \sqrt{2}\,\beta\,\bar{\textbf{q}}(t+t_\phi)^\intercal\check{\nabla}_{\hspace*{-0.04in}\perp}c}{\sqrt{\lvert \check{\nabla}_{\hspace*{-0.02in}\parallel}c\lvert^2+\beta^2\,\lvert \check{\nabla}_{\hspace*{-0.04in}\perp}c\lvert^2}},
\end{align}
where $\beta = \gamma \lvert\bm{\delta}(t)\lvert/(\mu\sqrt{2}\lvert \textbf{v}_m\lvert) = \gamma\omega_{\perp 0}/(\sqrt{2}\omega\mu\omega_{\parallel 0})$, $t_\phi = \phi/\omega$ with $\gamma$ and $\phi$ being gain and phase contribution of the linear part of the system (\ref{eq:pthway}) at the frequency $\omega = \lvert \bm{\Omega}_0\lvert$, $\check{\nabla} c = \nabla c(\bar{\textbf{p}})/\lvert \nabla c(\bar{\textbf{p}})\lvert$ is a unit vector in the direction of the gradient, and we used the following notations:
\begin{gather*}
\begin{aligned}
   \bar{\textbf{h}} &= \bar{\textbf{R}}\textbf{v}_m/\lvert \textbf{v}_m\lvert, & \bar{\textbf{q}}(t)&=  \bar{\textbf{R}}\bm{\delta}(t)/\lvert\bm{\delta}(t)\lvert,\\
   \check{\nabla}_{\hspace*{-0.02in}\parallel} c&= \bar{\textbf{h}}^\intercal \check{\nabla} c, & \check{\nabla}_{\hspace*{-0.04in}\perp} c&= \check{\nabla} c - \check{\nabla}_{\hspace*{-0.02in}\parallel} c\,\bar{\textbf{h}}.
\end{aligned}
\end{gather*}
Notably, the quasi-steady response is independent of the ambient concentration $c(\textbf{p}_0)$, the magnitude of the gradient $\lvert \nabla c(\bar{\textbf{p}})\lvert$, and the stimulus proportionality constant $\lambda$, all of which are irrelevant information from a chemotactic perspective. If we close the loop by replacing $\eta$ with the quasi-steady approximation $\eta_{\text{QS}}$, an intricate averaging analysis on the fast time scale $\tau = \omega t$ when $\omega \gg 1$ for the system of equations (\ref{eq:rot_kin})-(\ref{eq:pos_kin}) coupled with equation (\ref{eq:eta_qs}) leads to the following averaged quasi-steady equations:
\begin{subequations}\label{eq:qs_avg_dyn}
\begin{gather}\label{eq:qs_avg_dyn_1}
    \dot{\bar{\textbf{p}}}^\intercal\bar{\textbf{h}}=\frac{v\,\omega_{\parallel 0}}{\omega}\left(1+\frac{ \omega_{\perp 0}^2\omega_{\parallel 1}}{\alpha\,\omega^2\omega_{\parallel 0}}\check{\nabla}_{\hspace*{-0.02in}\parallel}c\right), \\\label{eq:qs_avg_dyn_2}
    \dot{\bar{\textbf{h}}}^\intercal\check{\nabla}c= \frac{\gamma\,\omega_{\perp 0}^2\omega_{\parallel 1}}{2\mu\alpha\omega^2\omega_{\parallel 0}}\cos(\phi) \lvert \check{\nabla}_{\hspace*{-0.04in}\perp}c\lvert^2,\\
    \alpha  = \sqrt{\lvert \check{\nabla}_{\hspace*{-0.02in}\parallel}c\lvert^2+\beta^2\lvert \check{\nabla}_{\hspace*{-0.04in}\perp}c\lvert^2}.
\end{gather}
\end{subequations}
Equation (\ref{eq:qs_avg_dyn_1}) expresses the speed along the average direction of motion $\bar{\textbf{h}}$, while equation (\ref{eq:qs_avg_dyn_2}) presents the rate of alignment of the average direction of motion $\bar{\textbf{h}}$ with the gradient.
We now analyze the qualitative dynamic behavior of the quasi-steady averaged equations (\ref{eq:qs_avg_dyn}) by considering three events and the corresponding response. The first event (the segments highlighted in green in Fig.\ref{fig:sample_traj}) is when the direction of average motion $\bar{\textbf{h}}$ is mostly aligned with the gradient (i.e. $\beta\lvert\check{\nabla}_{\hspace*{-0.04in}\perp} c\lvert \ll\check{\nabla}_{\hspace*{-0.02in}\parallel}c \approx 1$), in which case the response  $\eta_{\text{QS}}$ is approximately given by:
\begin{align}
     \eta_{\text{QS}}\approx 1 + \sqrt{2}\,\beta\,\bar{\textbf{q}}(t+t_\phi)^\intercal\check{\nabla}_{\hspace*{-0.04in}\perp}c,
\end{align}
where the second term is small compared to 1 (i.e. the periodic component is attenuated relative to the slope of the ramp component), and the change in the misalignment between the direction of average motion and the gradient is minor. Moreover, the average speed of the motion along the direction $\textbf{h}$ is increased:
\begin{align}\label{eq:speed_increase}
    \dot{\bar{\textbf{p}}}^\intercal \bar{\textbf{h}}
    &\approx \frac{v\,\omega_{\parallel 0}}{\omega}\left(1+\frac{ \omega_{\parallel 1}\omega_{\perp 0}^2}{\omega_{\parallel 0}\omega^2}\right),
\end{align}

The second event (the segments highlighted in purple in Fig.\ref{fig:sample_traj}) is when the direction of average motion is almost opposite to the gradient (i.e. $\check{\nabla}_{\hspace*{-0.02in}\parallel}c \approx -1$), in which case the response is approximately given by:
\begin{align}
     \eta_{\text{QS}}\approx -1 + \sqrt{2}\,\beta\,\bar{\textbf{q}}(t+t_\phi)^\intercal\check{\nabla}_{\hspace*{-0.04in}\perp}c,
\end{align} 
where once again the periodic term is small. However, the speed of the motion along the direction $\bar{\textbf{h}}$ is reduced:
\begin{align}\label{eq:speed_reduction}
    \dot{\bar{\textbf{p}}}^\intercal \bar{\textbf{h}}
    &\approx \frac{v\,\omega_{\parallel 0}}{\omega}\left(1-\frac{ \omega_{\parallel 1}\omega_{\perp 0}^2}{\omega_{\parallel 0}\omega^2}\right).
\end{align}
That is, when the motion is opposite to the gradient, the cell reduces its average speed along the direction of motion $\bar{\textbf{h}}$. This speed reduction mechanism can be observed in Fig.\ref{fig:sample_traj} as a gradual decrease in the helical pitch of the purple segment of the trajectory. Moreover, $\check{\nabla}_{\hspace*{-0.02in}\parallel}c \approx -1$ is an unstable direction for the average motion, so any slight misalignment triggers the transition towards the stable direction of motion  $\check{\nabla}_{\hspace*{-0.02in}\parallel}c \approx 1$. \\

The third event (the segments highlighted in red in Fig.\ref{fig:sample_traj}) is when the direction of average motion $\bar{\textbf{h}}$ is orthogonal to the gradient (i.e. $\check{\nabla}_{\hspace*{-0.02in}\parallel}c \approx 0$ and  $\lvert\check{\nabla}_{\hspace*{-0.04in}\perp}c\lvert \approx 1$), in which case the quasi-steady response $\eta_{\text{QS}}$ is dominated by the periodic component in the local concentration:
\begin{align}
    \eta_{\text{QS}}\approx  \sqrt{2}\,\bar{\textbf{q}}(t+t_\phi)^\intercal\check{\nabla}_{\hspace*{-0.04in}\perp}c,
\end{align}
and the alignment between the direction of average motion and the gradient is increased at a peak rate:
\begin{align} \label{eq:off_peak_rate}
    \dot{\bar{\textbf{h}}}^\intercal \check{\nabla} c &\approx \frac{\omega_{\perp 0}\omega_{\parallel 1}}{\sqrt{2}\,\omega}\cos(\phi).
\end{align}
We remark that near the maximum concentration, the gradient vanishes, and the behavior of the system is dominated by second order effects due to the Hessian of the concentration field which are neglected here. 
\section{Discussion}
Helical klinotaxis is a ubiquitous mode of taxis in microorganisms. In this study, we used sperm chemotaxis in sea urchins to highlight extremum seeking control as an underlying principle behind helical klinotaxis. This connection sheds light on the role played by the chemotactic signaling pathway and emphasizes the characterization of its dynamics as an adaptive band pass filter. Moreover, we showed that the switching-like behavior of sea urchin sperm \cite{jikeli2015sperm,kromer2018decision} can arise from a continuous dynamical description (\ref{eq:pthway}) without an explicit discontinuous switching logic as in previously proposed models. The key feature of the model (\ref{eq:pthway}) is that the adaptive gain $\rho$ adjusts according to the filtered stimulus $\dot{\xi}$ rather than the stimulus $s(t)$ directly. As a consequence, the ambient concentration levels do not alter the behavior of the model significantly. The forward speed of the cell is treated as a constant in the kinematic model (\ref{eq:orig_kin})-(\ref{eq:omega_def}). Yet, a cell is able to adjust its speed along the average direction of motion by dynamically regulating the angular velocity components. Our results suggest that this average speed reduction mechanism may be behind the peculiar switching-like behavior. That is, the on-response corresponds to the combined effect of speed increase and the attenuation of the periodic component when the direction of motion is almost parallel to the gradient. In contrast, the off-response may be explained as the combined effect of speed reduction when the direction of motion is opposite to the gradient followed by amplification of the periodic component when the direction of motion is misaligned with the gradient. The strength of the off-response is determined by the maximum speed reduction and the peak alignment rate given in equations (\ref{eq:speed_reduction}) and (\ref{eq:off_peak_rate}), respectively. In particular, the off response is most pronounced when $\omega_{\parallel 1}\omega_{\perp 0}^2 \approx \omega_{\parallel 0}\omega^2$, since it leads to a zero speed along the direction of motion $\bar{\textbf{h}}$ when it is opposite to the gradient. Furthermore, the feedback gain in the peak rate alignment depends on the factor $\cos(\phi)$, which attains its maximum value when the frequency of the periodic swimming pattern is inside the pass-band of the signaling pathway defined by $\mu $ and $\sigma$ so that the phase lag is minimal. Finally, we remark that the proposed connection between klinotaxis and extremum seeking may guide technological developments in robotic navigation \cite{long2004navigational,abdelgalil2022recursive}; it may inspire engineers to design source seeking algorithms with minimal sensors, suitable for miniaturized robots.}

\section*{Acknowledgements}
 MA conceptualized the work, performed theoretical analysis, simulations, and wrote a draft of the manuscript; YA and HT provided feedback on the physical conclusions and finalized the manuscript. HT likes to acknowledge the support of the NSF Grant CMMI-1846308.
\bibliography{sperm_chemotaxis}
\end{document}